\documentstyle[12pt]{article} 
\begin{document}
\newcommand{\singlespace}{
    \renewcommand{\baselinestretch}{1}
\large\normalsize}
\newcommand{\doublespace}{
   \renewcommand{\baselinestretch}{1.2}
   \large\normalsize}
\renewcommand{\theequation}{\thesection.\arabic{equation}}

\input amssym.def
\input amssym
\setcounter{equation}{0}
\def \ten#1{_{{}_{\scriptstyle#1}}}
\def \Z{\Bbb Z}
\def \C{\Bbb C}
\def \R{\Bbb R}
\def \Q{\Bbb Q}
\def \N{\Bbb N}
\def \l{\lambda}
\def \V{V^{\natural}}
\def \wt{{\rm wt}}
\def \tr{{\rm tr}}
\def \Res{{\rm Res}}
\def \End{{\rm End}}
\def \Aut{{\rm Aut}}
\def \mod{{\rm mod}}
\def \Hom{{\rm Hom}}
\def \im{{\rm im}}
\def \<{\langle} 
\def \>{\rangle} 
\def \w{\omega}
\def \c{{\tilde{c}}}
\def \o{\omega}
\def \t{\tau }
\def \ch{{\rm ch}}
\def \a{\alpha }
\def \b{\beta}
\def \e{\epsilon }
\def \la{\lambda }
\def \om{\omega }
\def \O{\Omega}
\def \qed{\mbox{ $\square$}}
\def \pf{\noindent {\bf Proof: \,}}
\def \voa{vertex operator algebra\ }
\def \voas{vertex operator algebras\ }
\def \p{\partial}
\def \1{{\bf 1}}
\def \ll{{\tilde{\lambda}}}
\def \H{{\bf H}}
\def \F{{\bf F}}
\def \h{{\frak h}}
\def \g{{\frak g}}
\def \rank{{\rm rank}}
\def \d{\delta}
\singlespace
\newtheorem{thmm}{Theorem}
\newtheorem{co}[thmm]{Corollary}
\newtheorem{th}{Theorem}
\newtheorem{prop}[th]{Proposition}
\newtheorem{coro}[th]{Corollary}
\newtheorem{lem}[th]{Lemma}
\newtheorem{rem}[th]{Remark}
\newtheorem{de}[th]{Definition}
\newtheorem{con}[th]{Conjecture}
\newtheorem{ex}[th]{Example}

\begin{center}
{\Large {\bf  Rational vertex operator algebras are finitely generated}} \\
\vspace{0.5cm}

Chongying Dong\footnote{Supported by NSF grants and a research grant from the
Committee on Research, UC Santa Cruz (dong@math.ucsc.edu).}  \ \ \  and \ \ \ 
Wei Zhang

\end{center}
\hspace{1.5 cm}

\begin{abstract} 
It is proved that any  vertex operator algebra 
for which  the image of the Virasoro element in Zhu's algebra
is algebraic over complex numbers is finitely generated. In particular,
any vertex operator algebra with a finite dimensional Zhu's algebra
is finitely generated. As a result, any rational vertex operator algebra is finitely generated.
\end{abstract}


Although many well known vertex operator algebras are finitely generated,
but whether or not an arbitrary  rational vertex operator algebra is finitely generated has been a basic problem in the theory of vertex operator algebra.
In this paper we give a positive answer to this problem and our result
justifies the assumption  in the physics literature that
any rational conformal field theory is finitely generated.

A systematic study of generators for an arbitrary vertex operator algebra was 
initiated in \cite{L}, \cite{KL}. A vertex operator algebra $V$ is called $C_1$-cofinite if $V=V_0+V_1\cdots+$ with $V_0$ being 1-dimensional
and $V/C_1(V)$ is finite dimensional where $C_1(V)$ is a subspace of $V$ 
spanned
by vectors $u_{-1}v, u_{-2}\1$ for $u,v\in V^+=\sum_{n>0}V_n$ and
$u_{n}$ is the component operator of $Y(u,z)=\sum_{n\in \Z}u_nz^{-n-1}.$ 
It is proved in \cite{L} that if a vertex operator algebra $V$ 
is  $C_1$-cofinite then it is finitely generated.  In fact,
if $\{x_i, i\in I\}$ is a set of vectors of $V$ such that $x_i+C_1(V)$ for $i\in I$
form a basis of $V/C_1(V),$ then $V$ is generated by this set of vectors
and $V$ has a PBW-like spanning set \cite{KL}. 

Another important finiteness for a vertex operator algebra is the
$C_2$-cofiniteness introduced by Zhu \cite{Z} in the proof of 
modular invariance of the $q$-characters of irreducible modules for 
a rational vertex operator algebra. A vertex operator algebra $V$ is called
$C_2$-cofinite if $V/C_2(V)$ is finite dimensional where
$C_2(V)$ is a subspace of $V$ spanned by $u_{-2}v$ for $u,v\in V.$ It is proved in \cite{Z} that
the span of the $q$-characters of irreducible modules for 
a rational, $C_2$-cofinite  vertex operator algebra affords a representation of the modular 
group $SL(2,\Z).$ So many results in the theory
of vertex operator algebra using the modularity of the $q$-characters of the
irreducible modules need both $C_2$-cofiniteness and rationality 
(see \cite{DLM3}, \cite{DM1}, \cite{DM2}, \cite{DM3}). It is shown in \cite{GN} that a $C_2$-cofinite
vertex operator algebra $V$ is finitely generated with a better PBW-like spanning set. Again one can choose a set $X$ of
homogeneous vectors of $V$  such that $x+C_2(V)$ for $x\in X$ 
form a basis of $V/C_2(V)$, then $V$ is spanned by
$$x^{1}_{-n_1}\cdots x^{k}_{-n_k}\1$$
where $x^i\in X$ and $n_1>n_2>\cdots >n_1>0.$ This result has been extended
in \cite{Bu} to give a spanning set for a weak module generated
by one vector.

The connection between the generators for a vertex operator algebra
and the $C_1$ or $C_2$-cofiniteness is not surprising as both $C_1$-cofinite
and $C_2$-cofinite properties are internal conditions on a vertex operator 
algebra. The mentioned results say that $C_1(V)$ or $C_2(V)$
can be generated by the other ``basic vectors.'' 
  
A vertex operator algebra is called {\em rational} if the admissible module
category is semisimple \cite{Z}, \cite{DLM2}, also see below.
It is clear that the rationality is an external condition. So the question
is how to connect the rationality and generators for a vertex operator 
algebra. The bridge is the associative algebra $A(V)$ introduced  in
\cite{Z}. We use the fact that the algebra $A(V)$ for rational vertex operator algebra $V$ 
is a finite dimensional semisimple associative algebra in this  paper 
to prove that any rational vertex operator algebra is finitely generated.

It turns out that a much weaker condition on $V$ is good enough to
guarantee that $V$ is finitely generated. Let $\omega$ be the Virasoro
vector of $V.$ Then the image $[\omega]$ of $\omega$ in $A(V)$ is a
central element. It is proved in this paper that if $[\omega]$ is
algebraic over $\C$ in $A(V)$ then $V$ is finitely generated. As
corollaries, rational vertex operator algebras are finitely generated
and their automorphism groups are algebraic groups \cite{DG}.

It is worthy to point out that although we can prove such vertex
operator algebra is  finitely generated, we do not known how to 
find a minimal set of generators for a rational vertex operator algebra.

We first review various notion of modules for a vertex operator 
algebra, following \cite{FLM}, \cite{Z} and \cite{DLM1}.

Let $V=(V,Y,\1,\omega)$ be a vertex operator algebra \cite{B}, \cite{FLM}.
A {\em weak} $V$ module is a vector space $M$ equipped with a linear map 
$$Y_M:V \rightarrow End(M)[[z,z^{-1}]]$$
$$v \mapsto Y_M(v,z)=\sum_{n \in \Z}v_n z^{-n-1}$$
and $v_n \in End(M)$. In 
addition $Y_M$ satisfies the following:

1) $v_nw=0$ for $n>>0$ where $v \in V$ and $w \in M$

2) $Y_M( {\textbf 1},z)=Id_M$

3) The Jacobi Identity holds:
\begin{eqnarray*}
& &z_0^{-1}\d ({z_1 - z_2 \over z_0})Y_M(u,z_1)Y_M(v,z_2)-
z_0^{-1} \d ({z_2- z_1 \over -z_0})Y_M(v,z_2)Y_M(u,z_1)\\
& &\ \ \ \ =z_2^{-1} \d ({z_1- z_0 \over z_2})Y_M(Y(u,z_0)v,z_2).
\end{eqnarray*}
\bigskip

An {\em admissible} $V$ module is a weak $V$ module which carries a
$\Z_+$ grading, $M=\bigoplus_{n \in \Z_+} M(n)$, such that if $v$ is homogeneousthen $v_m M(n) \subseteq M(n+\wt r-m-1).$ Since a uniform 
degree shift gives an isomorphic module we will assume $M(0)\ne 0$
if $M$ is nonzero.
\bigskip

An {\em ordinary} $V$ module is a weak $V$ module which carries a $\C$
grading, $M=\bigoplus_{\l \in \C} M_{\l}$, such that: 
1) $dim(M_{\l})< \infty,$ 2) $M_{\l+n=0}$ for fixed $\l$ and $n<<0,$
3) $L(0)w=\l w=wt(w) w$, for $w \in M$ where $L(0)$ is the component
operator of $Y_M(\omega,z)=\sum_{n\in\Z}L(n)z^{-n-2}.$

It is easy to see that an ordinary module is admissible. One of the main 
results in  \cite{DLM2} and \cite{Z} says that
if $V$ is rational then there are only finitely many irreducible
admissible modules up to isomorphism and each irreducible admissible module
is ordinary. 

A vertex operator algebra is called {\em regular} if every weak module
is a direct sum of irreducible ordinary modules (see \cite{DLM1}).  It
is clear that the regularity implies the rationality. One of the most
important problems in the theory of vertex operator algebra is to
understand the relations among rationality, regularity and
$C_2$-cofiniteness. It is shown that any regular vertex operator
algebra is $C_2$-cofinite and is finitely generated
\cite{L}. Conversely, a rational and $C_2$-cofinite vertex operator
algebra is regular \cite{ABD}. It is also known that a $C_2$-cofinite
vertex operator algebra is not necessarily rational \cite{A}. It is
suspected that rationality implies the $C_2$-cofiniteness but there is
no any progress in this direction so far.

For the purpose of the main result we also need to review the theory
of associative algebra $A(V)$ from \cite{Z}.

 For any homogeneous vectors $a\in V$, and $b\in V$, we define 
\begin{eqnarray*}
&a*b=\Res_{z}\frac{(1+z)^{\wt{a}}}{z}Y(a,z)b,\\ 
&a\circ b=\Res_{z}\frac{(1+z)^{\wt{a}}}{z^{2}}
Y(a,z)b,
\end{eqnarray*}
and extend to $V\times V$ bilinearly.  Denote by $O(V)$ the linear
span of $a\circ b$ ($a,b\in V$) and set $A(V)=V/O(V)$. Set $[a]=a+O(V)$ for
 $a\in V$.  Let $M$ be a weak $V$-module.
For a homogeneous $a\in V$ we write $o(a)$ for the operator $a_{\wt a-1}$
on $M.$  Extend the notation $o(a)$ to all $a\in V$ linearly. From the definition
of admissible module we see that $o(a)M(n)\subset M(n)$ for all $n\in \Z$ if
$M$ is an admissible module.  

 The following theorem is due to
\cite{Z} (also see \cite{DLM2}).
\begin{th}\label{P3.1}  Let $V$ be a vertex operator algebra. Then

(1) The bilinear operation $*$ induces $A(V)$ an associative algebra 
structure. The vector $[\1]$ is the identity and $[\w]$ is in the center of 
$A(V)$.

(2) Let $M=\bigoplus_{n=0}^{\infty}M(n)$ be an admissible $V$-module with
$M(0)\ne 0.$  
Then the linear map 
\[
o:V\rightarrow\End M(0),\;a\mapsto o(a)|_{M(0)}
\] 
induces an algebra homomorphism from $A(V)$ to $\End M(0)$. 
Thus $M(0)$ is a left $A(V)$-module. 

(3) The map $M\mapsto M(0)$ induces a bijection from the set of equivalence classes of irreducible admissible $V$-modules to the set of equivalence classes of irreducible $A(V)$-modules.    

(4) $[\omega]$ acts on any irreducible $A(V)$-module as a constant.

(5) If $V$ is rational then $A(V)$ is a finite dimensional semisimple
associative algebra.
\end{th}

A vector $u\in A(V)$ is called algebraic if there exits a nonzero
polynomial $f(x)\in \C[x]$ such that $f(u)=0.$ Note that if
$A(V)$ is finite dimensional then any $u\in A(V)$ is algebraic. 

We now in a position to prove that main result of this paper.
\begin{th}\label{mt} If $[\omega]\in A(V)$ is algebraic then $V$ is finitely
generated. In particular, if $A(V)$ is finite dimensional, $V$ is finitely generated. Furthermore, any rational vertex operator algebra
is finitely generated.  
\end{th}

\pf Since $[\omega]\in A(V)$ is algebraic then $[\omega]$ satisfies
$$([\omega]-\l_1)\cdots ([\omega]-\l_k)=0$$
for some complex numbers $\l_i.$ That is, there exist
homogeneous vectors $u^i,v^i\in V$ for $i=1,...,p$ such that
$$(\omega-\l_1)*\cdots *(\omega-\l_k)=\sum_{i}u^i\circ v^i.$$
Let $n$ be a positive integer such that $n>\l_j$ if
$\l_j\in \Z$ and that $n$ is greater than or equal to the weights
of $u^i,v^i$ for all $i.$ We can also assume that $n\geq 2.$ 
Let $U$ be a vertex operator subalgebra generated by $\sum_{m\leq n}V_m.$
We claim that $V=U.$ 

If $U\ne V,$ let $t$ be the minimal positive integer such that
$V_t\ne U_t.$ Then $V/U=\oplus_{m\geq t}V_m/U_m$ is a $U$-module 
where $U_m=V_m\cap U.$ Note that $V_t/U_t$ is an $A(U)$-module and
$[\omega]_U=\omega+O(U)$ acts on $V_t/U_t$ as $t.$ On the other hand,
$[\omega]_U$ is also algebraic in $A(U)$ and satifies the same
relation
$$([\omega]_U-\l_1)\cdots ([\omega]_U-\l_k)=0.$$ 
That is, if $[\omega]_U$ acts on an $A(U)$-module as a constant,
this constant must be one of $\l_j.$ Since $t$ is different from
$\l_j$ for all $j$ we have a contradiction. \qed
 
It is worthy to note that if $V$ is $C_2$-cofinite, then $A(V)$ is finite
dimensional. In fact, it is easy to verify that if $x_i\in V$ $(i\in I)$ be a set of
homogeneous vectors such that $x_i+C_2(V)$ for $i\in I$ form a basis
of $V/C_2(V)$ then $x_i+O(V)$ for $i\in I$ give a spanning set
of $A(V).$ In particular, $\dim A(V)\leq \dim V/C_2(V).$  So an immediate
corollary is that a $C_2$-cofinite vertex operator algebra is finitely 
generated. 

We can strengthen Theorem \ref{mt} in the following way.

\begin{coro} Let $V$ be a vertex operator algebra and $U$ a vertex operator
subalgebra with the same Virasoro algebra. If $\omega+O(U)$ is algebraic
in $A(U)$ then $V$ is finitely generated. In particular, any vertex operator
algebra which has a rational vertex operator subalgebra is finitely generated.
\end{coro}

\pf It is clear that $O(U)$ is a subspace of $O(V).$ So the embedding 
from $U$ to $V$ induces an algebra homomorphism. So $[\omega]$ is also 
algebraic in $A(V).$ By Theorem \ref{mt} we have the result. \qed

Another application of Theorem \ref{mt} is on the automorphism group of
a rational vertex operator algebra. 
Recall that an automorphism $g$ of a vertex operator algebra $V$
is a linear isomorphism from $V$ to $V$ such that $g\1=\1,$ 
$g\omega=\omega$ and $gY(v,z)g^{-1}=Y(gv,z)$ for all $v\in V.$ Let
$\Aut(V)$ denote the full automorphism group.

\begin{coro} If $A(V)$ is finite dimensional then $\Aut(V)$ is an algebraic
group. In particular, the automorphism group of any rational
vertex operator algebra is an algebraic group.
\end{coro}

\pf The corollary follows from Theorem \ref{mt} and a result in \cite{DG}
which says the automorphism group of 
any finitely generated vertex operator algebra is an algebraic group.
\qed

\end{document}